\documentclass[12pt]{amsart}

\newtheorem{lemma}{Lemma}

\newtheorem*{theorem}{Theorem}
\newtheorem*{corollary}{Corollary}

\newcommand{\bC}{\mathbb{C}}
\newcommand{\bQ}{\mathbb{Q}}
\newcommand{\bR}{\mathbb{R}}
\newcommand{\bZ}{\mathbb{Z}}
\newcommand{\cC}{{\mathcal{C}}}
\newcommand{\Res}{{\mathrm{\,Res}}}

\begin{document}

\baselineskip=17pt

\title{A Gel'fond type criterion in degree two}
\author{Benoit Arbour}
\address{
   Department of Mathematics\\
   McGill University\\
   805 Sherbrooke Ouest\\
   Montr\'eal, Qu\'ebec H3A 2K6, Canada}
\email{arbour@math.mcgill.ca}
\author{Damien ROY}
\address{
   D\'epartement de Math\'ematiques\\
   Universit\'e d'Ottawa\\
   585 King Edward\\
   Ottawa, Ontario K1N 6N5, Canada}
\email{droy@uottawa.ca}
\thanks{Work partially supported by NSERC and CICMA}

\maketitle


\section{Introduction}


Let $\xi$ be any real number and let $n$ be a positive integer.
Defining the {\it height} $H(P)$ of a polynomial $P$ as the
largest absolute value of its coefficients, an application of
Dirichlet box principle shows that, for any real number $X\ge 1$,
there exists a non-zero polynomial $P\in\bZ[T]$ of degree at most
$n$ and height at most $X$ which satisfies
$$
|P(\xi)| \le c X^{-n}
$$
for some suitable constant $c>0$ depending only on $\xi$ and $n$.
Conversely, Gel'fond's criterion implies that there are constants
$\tau=\tau(n)$ and $c=c(\xi,n)>0$ with the property that if, for
any real number $X\ge 1$, there exists a non-zero polynomial
$P\in\bZ[T]$ with
$$
 \deg(P) \le n,
   \quad
 H(P)\le X
   \quad \hbox{and} \quad
 |P(\xi)|\le cX^{-\tau}
$$
then $\xi$ is algebraic over $\bQ$ of degree at most $n$.  For
example, Brownawell's version of Gel'fond's criterion in
\cite{Bro1} implies that the above statement holds with any
$\tau>3n$, and the more specific version proved by Davenport and
Schmidt as Theorem 2b of \cite{DS2} shows that it holds with
$\tau=2n-1$. On the other hand, the above application of Dirichlet
box principle implies $\tau \ge n$. So, if we denote by $\tau_n$
the infimum of all admissible values of $\tau$ for a fixed $n\ge
1$, then we have $\tau_1=1$ and, in general
$$
n \le \tau_n \le 2n-1.
$$

In the case of degree $n=2$, the study of a specific class of
transcendental real numbers in \cite{Roy1} provides the sharper
lower bound $\tau_2\ge \gamma^2$ where $\gamma =(1 + \sqrt{5})/2$
denotes the golden ratio (see Theorem 1.2 of \cite{Roy1}). Our
main result below shows that we in fact have $\tau_2=\gamma^2$ by
establishing the reverse inequality $\tau_2\le \gamma^2$:

\begin{theorem}

Let $\xi \in \bC$. Assume that for any sufficiently large positive
real number $X$ there exists a non-zero polynomial $P \in \bZ[T]$
of degree at most $2$ and height at most $X$ such that
    \begin{equation}\label{hypothesis}
        |P(\xi)| \leq \frac{1}{4}X^{-\gamma^2}.
    \end{equation}
Then $\xi$ is algebraic over $\bQ$ of degree at most $2$.

\end{theorem}

Comparing this statement with Theorem 1.2 of \cite{Roy1}, we see
that it is optimal up to the value of the multiplicative constant
$1/4$ in (\ref{hypothesis}).  Although we do not know the best
possible value for this constant, our argument will show that it
can be replaced by any real number $c$ with $0 < c < c_0=(6 \cdot
2^{\frac{1}{\gamma}})^{-\frac{1}{\gamma}} \cong 0.253$. As the
reader will note, our proof, given in Section 3 below, has the
same general structure as the proof of the main result of
\cite{DS1} and the proof of Theorem 1a of \cite{DS2}.

\medskip
Following the method of Davenport and Schmidt in \cite{DS2}
combined with ideas from \cite{BugTeu1} and \cite{RoyWal1}, we
deduce the following result of simultaneous approximation of a
real number by conjugate algebraic numbers:

\begin{corollary}

Let $\xi$ be a real number which is not algebraic over $\bQ$ of
degree at most $2$. Then, there are arbitrarily large real numbers
$Y\ge 1$ for which there exist an irreducible monic polynomial $P
\in \bZ[T]$ of degree $3$ and an irreducible polynomial $Q \in
\bZ[T]$ of degree 2, both of which have height at most $Y$ and
admit at least two distinct real roots whose distance to $\xi$ is
at most $cY^{-(3-\gamma)/2}$, with a constant $c$ depending only
on $\xi$.

\end{corollary}

The proof of this corollary is postponed to Section 4 below.


\section{Preliminaries}


We collect here several lemmas which we will need in the proof of
the theorem. The first one is a special case of the well-known
Gel'fond's lemma for which we computed the optimal values of the
constants.

\begin{lemma}\label{gelfond}
Let $L, M \in \bC[T]$ be polynomials of degree at most $1$. Then
we have
    \begin{equation*}
        \frac{1}{\gamma} H(L)H(M) \leq H(LM) \leq 2H(L)H(M).
    \end{equation*}
\end{lemma}

The second result is an estimate for the resultant of two
polynomials of small degree.

\begin{lemma}\label{two:resultant}
Let $m,n \in \{1,2\}$, and let $P$ and $Q$ be non-zero polynomials
in $\bZ[T]$ with $\deg(P)\leq m$ and $\deg(Q)\leq n$. Then, for
any complex number $\xi$, we have
    \begin{equation*}
        |\Res(P,Q)| \leq H(P)^n H(Q)^m
                        \Big( c(m,n) \frac{|P(\xi)|}{H(P)}
                            + c(n,m) \frac{|Q(\xi)|}{H(Q)} \Big)
    \end{equation*}
where $c(1,1)=1$, $c(1,2)=3$, $c(2,1)=1$ and $c(2,2)=6$.
\end{lemma}

The proof of the above statement is easily reduced to the case
where $\deg(P)=m$ and $\deg(Q)=n$.  The conclusion then follows by
writing $\Res(P,Q)$ as a Sylvester's determinant and by arguing as
Brownawell in the proof of lemma 1 of \cite{Bro1} to estimate this
determinant.

\medskip
The third lemma may be viewed, for example, as a special case of
Lemma 13 of \cite{LR}.

\begin{lemma}\label{two:root}
Let $P,Q \in \bZ[T]$ be non-zero polynomials of degree at most $2$
with greatest common divisor $L\in\bZ[T]$ of degree $1$. Then, for
any complex number $\xi$, we have
  \begin{equation*}
        H(L)|L(\xi)|
        \leq \gamma \big( H(P)|Q(\xi)| + H(Q)|P(\xi)| \big).
  \end{equation*}
\end{lemma}

\begin{proof}
The quotients $P/L$ and $Q/L$ being relatively prime polynomials
of $\bZ[T]$, their resultant is a non-zero integer. Applying Lemma
\ref{two:resultant} with $m=n=1$ and using Lemma \ref{gelfond}, we
then deduce, if $L(\xi)\neq 0$,
  \begin{eqnarray}
       1 \leq |\Res(P/L,Q/L)|
         &\leq& H (P/L) \big| (Q/L)(\xi) \big|
              + H (Q/L)\big| ( P/L)(\xi) \big| \nonumber \\
         &\leq& \gamma \frac{H(P)}{H(L)} \frac{|Q(\xi)|}{|L(\xi)|}
              + \gamma \frac{H(Q)}{H(L)} \frac{|P(\xi)|}{|L(\xi)|}\ .
              \nonumber
  \end{eqnarray}
\end{proof}

\begin{lemma}\label{three:sum}
Let $\xi\in\bC$ and let $P,Q,R \in \bC[T]$ be arbitrary
polynomials of degree at most $2$. Then, writing the coefficients
of these polynomials as rows of a $3\times 3$ matrix, we have
  \begin{equation*}
    |\det(P,Q,R)| \leq 2H(P)H(Q)H(R) \biggl( \frac{|P(\xi)|}{H(P)}
                      + \frac{|Q(\xi)|}{H(Q)}
                      + \frac{|R(\xi)|}{H(R)}\biggr)
  \end{equation*}
\end{lemma}

This above statement follows simply by observing, as in the proof
of Lemma 4 of \cite{DS1}, that the determinant of the matrix does
not change if, in this matrix, we replace the constant
coefficients of $P$, $Q$ and $R$ by the values of these
polynomials at $\xi$.

\medskip
We also construct a sequence of ``minimal polynomials'' similarly
as it is done in \S3 of \cite{DS1}:

\begin{lemma}\label{definition:P}
Let $\xi \in \bC$ with $[\bQ(\xi):\bQ] >2$. Then there exists a
strictly increasing sequence of positive integers $(X_i)_{i\geq1}$
and a sequence of non-zero polynomials $(P_i)_{i\geq1}$ in
$\bZ[T]$ of degree at most $2$ such that, for each $i\geq 1$, we
have:
  \begin{itemize}
    \item $H(P_i) = X_i$,
    \item $|P_{i+1}(\xi)| < |P_i(\xi)|$,
    \item $|P_i(\xi)|\leq |P(\xi)|$ for all $P \in \bZ[T]$
            with $deg(P) \leq 2$ and $0<H(P)<X_{i+1}$,
    \item $P_i$ and $P_{i+1}$ are linearly independent
            over $\bQ$.
  \end{itemize}
\end{lemma}

\begin{proof}
For each positive integer $X$, define $p_X$ to be the smallest
value of $|P(\xi)|$ where $P \in \bZ[T]$ is a non-zero polynomial
of degree $\leq 2$ and height $\leq X$. This defines a
non-decreasing sequence $p_1 \geq p_2 \geq p_3 \geq \dots$ of
positive real numbers converging to $0$. Consider the sequence
$X_1 < X_2 < \dots$ of indices $X\ge 2$ for which $p_{X-1} >
p_{X}$. For each $i \geq 1$, there exists a polynomial $P_i \in
\bZ[T]$ of degree $\leq 2$ and height $X_i$ with $|P_i(\xi)|=p_i$.
The sequences $(X_i)_{i\geq1}$ and $(P_i)_{i\geq1}$ clearly
satisfy the first three conditions. The last condition follows
from the fact that the polynomials $P_i$ are primitive of distinct
height.
\end{proof}

\begin{lemma}\label{three:independent}
Assume, in the notation of Lemma \ref{definition:P}, that
  \begin{equation*}
    \lim_{i \to \infty} X_{i+1}|P_i(\xi)| = 0.
  \end{equation*}
Then there exist infinitely many indices $i\ge 2$ for which
$P_{i-1}$, $P_i$ and $P_{i+1}$ are linearly independent over
$\bQ$.
\end{lemma}

\begin{proof}
Assume on the contrary that $P_{i-1}, P_i$ and $ P_{i+1}$ are
linearly dependent over $\bQ$ for all $i \geq i_0$. Then the
subspace $V$ of $\bQ[T]$ generated by $P_{i-1}$ and $P_i$ is
independent of $i$ for $i \geq i_0$.  Let $\{P,Q\}$ be a basis of
$V \cap \bZ^3$. Then, for each $i\ge i_0$, we can write
  \begin{equation*}
        P_i =a_iP + b_iQ
  \end{equation*}
for some integers $a_i$ and $b_i$ of absolute value at most $c
X_i$, with a constant $c>0$ depending only on $P$ and $Q$.  Since
$P_i$ and $P_{i+1}$ are linearly independent, we get
  \begin{equation*}
    1 \leq \left\| \begin{array}{cc}
            a_i & b_i \\
            a_{i+1} & b_{i+1}
        \end{array} \right\|
    = \frac{| a_i P_{i+1}(\xi) - a_{i+1} P_i(\xi) |}{|Q(\xi)|}
    \leq \frac{2c}{|Q(\xi)|}X_{i+1}|P_i(\xi)|
    \end{equation*}
in contradiction with the hypothesis as we let $i$ tend to
infinity.
\end{proof}


\section{Proof of the theorem}


Let $c$ be a positive real number and let $\xi$ be a complex
number with $[\bQ(\xi):\bQ]>2$.  Assume that, for any sufficiently
large real number $X$, there exists a non-zero polynomial $P \in
\bZ[T]$ of degree $\leq 2$ and height $\leq X$ with $|P(\xi)| \leq
cX^{-\gamma^2}$.  We will show that these conditions imply $c \ge
c_0 = (6\cdot 2^{1/\gamma})^{-1/\gamma} >1/4$, thereby proving the
theorem.

Let $c_1$ be an arbitrary real number with $c_1>c$.  By virtue of
our hypotheses, the sequences $(X_i)_{i \geq 1}$ and $(P_i)_{i
\geq 1}$ given by Lemma \ref{definition:P} satisfy
    \begin{equation*}
        |P_i(\xi)|\leq cX_{i+1}^{-\gamma^2}
    \end{equation*}
for any sufficiently large $i$. Then, by Lemma
\ref{three:independent}, there exists infinitely many $i$ such
that $P_{i-1}$, $P_i$ and $P_{i+1}$ are linearly independent. For
such an index $i$, the determinant of these three polynomials is a
non-zero integer and, applying Lemma \ref{three:sum}, we deduce
    \begin{eqnarray}
        1 \leq |\det(P_{i-1},P_i,P_{i+1})|
        &\leq& 2X_{i-1}X_iX_{i+1}
           \biggl( \frac{|P_{i-1}(\xi)|}{X_{i-1}}
                 + \frac{|P_{i}(\xi)|}{X_{i}}
                 + \frac{|P_{i+1}(\xi)|}{X_{i+1}} \biggr) \nonumber \\
        &\leq&2cX_i^{-\gamma}X_{i+1} +4cX_{i+1}^{1-\gamma}. \nonumber
    \end{eqnarray}
Assuming that $i$ is sufficiently large, this implies
    \begin{equation} \label{first:estimate}
    X_i^{\gamma} \leq 2c_1 X_{i+1}.
    \end{equation}

Suppose first that $P_i$ and $P_{i+1}$ are not relatively prime.
Then, their greatest common divisor is an irreducible polynomial
$L\in\bZ[T]$ of degree $1$, and Lemma \ref{two:root} gives
    \begin{equation} \label{second:estimate}
    H(L)|L(\xi)|
      \leq \gamma \Big(X_i |P_{i+1}(\xi)| + X_{i+1}|P_i(\xi)|\Big)
      \leq 2 \gamma c X_{i+1}^{-\gamma}.
    \end{equation}
Since $P_{i-1}$, $P_i$ and $P_{i+1}$ are linearly independent, the
polynomial $L$ does not divide $P_{i-1}$ and so the resultant of
$P_{i-1}$ and $L$ is a non-zero integer.  Applying Lemma
\ref{two:resultant} then gives
    \begin{eqnarray}
        1 \leq |\Res(P_{i-1},L)|
        &\leq& H(P_{i-1}) H(L)^2
               \left( \frac{|P_{i-1}(\xi)|}{H(P_{i-1})}
                 + 3 \frac{|L(\xi)|}{H(L)} \right)
               \nonumber \\
        &\leq& cX_i^{-\gamma^2}H(L)^2 + 3X_{i-1} H(L) |L(\xi)|.
               \nonumber
    \end{eqnarray}
Combining this with (\ref{second:estimate}) and with the estimate
$H(L)\le \gamma H(P_i) \le \gamma X_i$ coming from Lemma
\ref{gelfond}, we conclude that, in this case, the index $i$ is
bounded.

Thus, assuming that $i$ is sufficiently large, the polynomials
$P_i$ and $P_{i+1}$ are relatively prime and therefore their
resultant is a non-zero integer.  Using Lemma \ref{two:resultant}
we then find
    \begin{equation*}
    1 \leq |\Res(P_i,P_{i+1})|
    \leq 6 X_iX_{i+1}
         \big(cX_iX_{i+2}^{-\gamma^2} + cX_{i+1}^{-\gamma}\big)
    \leq 6c_1X_iX_{i+1}^{1-\gamma}
    \end{equation*}
since from (\ref{first:estimate}), we have $cX_i \leq
(c_1-c)X_{i+1}$ for large $i$. Using (\ref{first:estimate}) again,
this implies
    \begin{equation*}
        1 \leq 6c_1(2c_1)^{1/\gamma}
    \end{equation*}
and thus, $c_1 \geq c_0=(6 \cdot 2^{1/\gamma})^{-1/\gamma}$. The
choice of $c_1>c$ being arbitrary, this shows that $c\ge c_0$ as
announced.


\section{Proof of the corollary}


Let $\xi$ be as in the statement of the corollary and let $V$
denote the real vector space of polynomials of degree at most $2$
in $\bR[T]$.  It follows from the theorem that there exist
arbitrarily
large real numbers $X$ for which the convex body
$\cC(X)$ of $V$
defined
by
  \begin{equation*}
  \cC(X) = \{ P\in V \,;\, |P(\xi)|\le (1/4)X^{-\gamma^2},\
           |P'(\xi)| \le c_1X \ \mbox{and}\ |P''(\xi)|\le c_1X \}
  \end{equation*}
with $c_1=(1+|\xi|)^{-2}$ contains no non-zero integral
polynomial.  By Proposition 3.5 of \cite{RoyWal1} (a version of
Mahler's theorem on polar reciprocal bodies), this implies that
there exists a constant $c_2> 1$ such that, for the same values of
$X$, the convex body
  \begin{equation*}
  \cC^*(X) = \{ P\in V \,;\, |P(\xi)|\le c_2X^{-1},\
               |P'(\xi)| \le c_2X^{-1} \ \mbox{and}\
               |P''(\xi)|\le c_2X^{\gamma^2} \}
  \end{equation*}
contains a basis of the lattice of integral polynomials in $V$.

Fix such an $X$ with $X\ge 1$, and let $\{P_1,P_2,P_3\}\subset
\cC^*(X)$ be a basis of $V\cap\bZ[T]$.  We now argue as in the
proof of Proposition 9.1 of \cite{RoyWal1}.  We put
  \begin{equation*}
    B(T)=T^2-1, \quad
    r=X^{-(1+\gamma^2)/2} \quad\mbox{and}\quad
    s=20c_2X^{-1},
  \end{equation*}
and observe that any polynomial $S\in V$ with $H(S-B)< 1/3$ admits
at least two real roots in the interval $[-2,2]$ as such a
polynomial takes positive values at $\pm 2$ and a negative value
at $0$. We also note that, since $P_i\in\cC^*(X)$, we have
  \begin{equation*}
  H(P_i(rT+\xi)) \le c_2X^{-1},\quad (i=1,2,3).
  \end{equation*}
Since $\{P_1,P_2,P_3\}$ is a basis of $V$ over $\bR$, we may write
  \begin{equation*}
    (T-\xi)^3 + sB\Big(\frac{T-\xi}{r}\Big)
      = T^3 + \sum_{i=1}^3 \theta_i P_i(T)
    \quad\hbox{and}\quad
    sB\Big(\frac{T-\xi}{r}\Big)
      = \sum_{i=1}^3 \eta_i P_i(T)
  \end{equation*}
for some real numbers $\theta_1,\theta_2,\theta_3$ and
$\eta_1,\eta_2,\eta_3$.  For $i=1,2,3$, choose integers $a_i$ and
$b_i$ with $|a_i-\theta_i|\le 2$ and $|b_i-\eta_i|\le 2$ so that
the polynomials
  \begin{equation*}
    P(T) = T^3 + \sum_{i=1}^3 a_iP_i(T)
    \quad\hbox{and}\quad
    Q(T) = \sum_{i=1}^3 b_iP_i(T)
  \end{equation*}
are respectively congruent to $T^3+2$ and $T^2+2$ modulo $4$.
Then, by Eisenstein's criterion, $P$ and $Q$ are irreducible
polynomials of $\bZ[T]$.  Moreover, we find
  \begin{eqnarray}
  H\big(s^{-1}P(rT+\xi)-B(T)\big)
    &=& s^{-1}
      H\Big( (rT)^3 + \sum_{i=1}^3 (a_i-\theta_i)P_i(rT+\xi)
       \Big) \nonumber \\
    &\le& s^{-1}\max\{r^3, 6c_2X^{-1}\} \nonumber \\
    &<& 1/3. \nonumber
  \end{eqnarray}
Then, $P(rT+\xi)$ has at least two distinct real roots in the
interval $[-2,2]$ and so $P$ has at least two real roots whose
distance to $\xi$ are at most $2r$. A similar but simpler
computation shows that the same is true of the polynomial $Q$.
Finally, the above estimate implies $H(P(rT+\xi)) \le 4s/3$ and so
$H(P) \le c_3X^{\gamma^2}$ for some constant $c_3>0$, and the same
for $Q$.  These polynomials thus satisfy the conclusion of the
corollary with $Y=c_3X^{\gamma^2}$ and an appropriate choice of
$c$.



\end{document}